\magnification=\magstep 1
  \input amstex
\input epsf

\documentstyle{onecolumn}
%

\ifx\macrosamsLoaded\relax 
\else\global\let\macrosamsLoaded\relax\fi

    \ifx\everymathLoaded\relax
    \else

\ifx\everymathLoaded\relax 
\else\global\let\everymathLoaded\relax\fi

\def\idtb{\mathrel{:\mkern-2.0mu=}}  

\def\lcm{\operatorname{lcm}}

\def\ut{{\mkern-1mu\times\mkern-2.5mu}}		
\def\ct{{\mkern-1.7mu\times\mkern-1.9mu}}

\def\0{{\bold0}}	\def\1{{\bold1}}
\let\dirac=\1  
\def\ind#1{{\bold1}_{#1}}	


\def\negskip{\mkern -1mu}
\xdef\doubledots{\nobreak\, .\nobreak\negskip. \nobreak\,}
\def\ico#1,#2.{ [ #1 \doubledots #2 ) }
\def\icc#1,#2.{ [ #1 \doubledots #2 ] }
\def\ioc#1,#2.{ ( #1 \doubledots #2 ] }
\def\ioo#1,#2.{ ( #1 \doubledots #2 ) }

\def\Ico#1,#2.{ \bigl[ #1 \doubledots #2 \bigr) }
\def\Icc#1,#2.{ \bigl[ #1 \doubledots #2 \bigr] }
\def\Ioc#1,#2.{ \bigl( #1 \doubledots #2 \bigr] }
\def\Ioo#1,#2.{ \bigl( #1 \doubledots #2 \bigr) }

\def\SOI#1,#2.{ #1 \doubledots #2 }	

\def\intervalspacing{,}
\def\rco#1,#2.{ [ #1 \intervalspacing #2 ) }
\def\rcc#1,#2.{ [ #1 \intervalspacing #2 ] }
\def\roc#1,#2.{ ( #1 \intervalspacing #2 ] }
\def\roo#1,#2.{ ( #1 \intervalspacing #2 ) }

\def\Rco#1,#2.{ \bigl[ #1 \intervalspacing #2 \bigr) }
\def\Rcc#1,#2.{ \bigl[ #1 \intervalspacing #2 \bigr] }
\def\Roc#1,#2.{ \bigl( #1 \intervalspacing #2 \bigr] }
\def\Roo#1,#2.{ \bigl( #1 \intervalspacing #2 \bigr) }


\define\ii{\infty}
%


\redefine\hom#1{{  \operatorname{Hom}(#1)  }} 
	\let\emptyset=\varnothing
\let\epsilon=\varepsilon    
\let\temp=\phi      \let\phi=\varphi          \let\varphi=\temp

\def\complex{{\Bbb C}}
\def\integers{  {\Bbb Z}  }

\def\rationals{  {\Bbb Q}  }






%

%

\define\setdiff{\smallsetminus}

			
\let\la=\langle			\let\ra=\rangle
\fi

\catcode`\@=11

%
%
%



\define\scap{\raise6.75pt\hbox{\vrule height.85pt depth8.75pt width.95pt\vrule
   width 5.1pt height.95pt depth0pt\hfill\vrule height1pt depth8.5pt width.95pt
}}
\define\bscap{\raise8.5pt\hbox{\vrule height.85pt depth13.75pt width.95pt\vrule
   width 8pt height.95pt depth0pt\hfill\vrule height.9pt depth13.5pt width.95pt
}}

\newif\ifdraft	\draftfalse	




\define\restrict#1{|_{#1}}
\define\ondots#1{\overset{#1}\to\ldots}
\define\oncdots#1{\overset{#1}\to\cdots}

\define\udots{\mathinner{\lower0\p@\vbox{\kern7\p@\hbox{.}}\mkern2mu
    \lower-3\p@\hbox{.}\mkern2mu\lower-6\p@\hbox{.}\mkern1mu}}

\define\elementdiff{\mathbin{\hbox{$\mkern4mu\raise0.79ex\hbox{$\scriptscriptstyle \in$}\mkern-10mu\hbox{$-$}$}}}



\define\conj#1{{#1}^*}		
\define\dual#1{{#1}^{\mkern-3mu\divideontimes}}		

\define\eqnote{{ \mathrel{\;\overset\text{note}\to{=\mkern-3.1mu =\mkern-3.1mu =}\;}  }}


\define\smallpmatrix#1\endsmallpmatrix{%
  \left(\smallmatrix#1\endsmallmatrix\right)%
  }
\define\smallbmatrix#1\endsmallbmatrix{%
  \left[\smallmatrix#1\endsmallmatrix\right]%
  }


\newcount\internalcounter

\def\initctr#1 {\internalcounter=#1}	\initctr1	
\def\putctr#1{\relax
  \xdef#1{\number\internalcounter}\global\advance\internalcounter by1
  \xdef\currentnumber{#1}%
  #1{\ifdraft\text{\eighttt\string#1}\fi}}
\def\puteqn#1{%
  \global\xdef#1{\huchapternumber.\number\internalcounter}%
  \global\advance\internalcounter by1
  \xdef\currentnumber{#1}%
  #1{\ifdraft\text{\tt\string#1}\fi}}
\let\puttype=\putctr 
\let\put=\putctr 

\newwrite\pfchannel		
\xdef\currentfile{\jobname}	
\def\inputsource#1 {
  \closeout\pfchannel	
  \xdef\currentfile{#1}%
  \initializepf		
  \input\currentfile
  }
\def\initializepf{
  \global\def\pf{%
    \immediate\write16{Opening forward reference: refs.\currentfile.tex}
    \immediate\openout\pfchannel=refs.\currentfile.tex \unskip
    \global\let\pf=\internalpf%
    \pf%
    }%
  }
\initializepf	

\def\internalpf#1{%
  \puttype#1
  \immediate\write\pfchannel{\def\noexpand#1{#1}}%
  \immediate\write16{\def\noexpand#1{#1}}%
  }


\def\display#1\enddisplay{\dimen1=.9\hsize
  \setbox0=\hbox{\proclaimfont\ignorespaces#1\unskip}	\dimen2=\wd0 
  \advance\dimen2 by \dimen1		\divide\dimen2 by \dimen1
    \count255=\dimen2
  \dimen2=\wd0				\divide\dimen2 by \count255
  \foldedtext\foldedwidth{\dimen2}{\unhbox0}
  }

\long\def\chapterhead#1#2\endchapterhead{
  \messageln{CHAPTER #1 #2}\goodbreak\bigskip 
  \vskip0pt plus 30pt \goodbreak
  \xdef\huchapternumber{#1}%
  \xdef\huchaptertext{#2}%
  \centerline{\twelvepoint\sc \S #1 \ \ \ignorespaces #2}%
  \nobreak\medskip\noindent\tenpoint\ignorespaces
  }

\def\huchapternumber{0}	
\long\def\huchapter#1#2{\messageln{CHAPTER #1 #2}\goodbreak\bigskip 
  \vskip0pt plus 30pt \goodbreak 
  \xdef\huchapternumber{#1}\initctr1
  \xdef\huchaptertext{#2}
  \centerline{\twelvepoint\sc \S #1 \ \ \ignorespaces #2}%
  \nobreak\medskip\noindent\tenpoint\ignorespaces
  }


\def\specialheadheight{1.5pt}

\def\specialhead#1\endspecialhead{%
\goodbreak
\bigskip\begingroup 
  \baselineskip=0pt
  \lineskiplimit=0pt
  \lineskip=4pt
  \setbox0\hbox{\ignorespaces\bf#1\unskip}
  \setbox1\hbox to \wd0{\leaders\hrule
    height \specialheadheight \hfill}
  \line{\hfil \vbox{\copy1\copy0\copy1}}
  \nobreak\vskip 2.0pt minus 1.0pt
\endgroup\noindent
}

\def\largehead#1\endlargehead{\goodbreak\bigskip
  \centerline{\twelvepoint\sc\ignorespaces#1}%
  \nobreak\medskip\noindent\tenpoint\ignorespaces
  }

\def\variablehead#1\endvariablehead{\goodbreak\bigskip
  \centerline{\twelvepoint\sc\ignorespaces#1}%
  \nobreak\medskip\noindent\tenpoint\ignorespaces
  }

\comment
\def\section#1{
  \medbreak\smallskip
  \noindent{\bf #1}\hskip12pt plus6pt minus4pt\ignorespaces
  }
\def\subsection#1{\smallbreak\smallskip
  {\bit #1}\hskip8pt plus4pt minus4pt\ignorespaces 
  }
\endcomment

\def\startnormal{\par
  \medskip
  \noindent \rm
  }

\def\?#1{\text{\itemfont#1}}		\def\=#1{({\itemfont#1\/})}

\def\proof{\medbreak\smallskip
  \noindent\kern12pt{\demofont@ Proof.}\hskip8pt plus6pt minus4pt\ignorespaces
  }

\def\proofof#1{\medbreak\smallskip
  \noindent\kern12pt{\demofont@ Proof of ({\itemfont#1}).}\hskip8pt plus6pt minus4pt\ignorespaces
  }

\def\proofwith#1{\medbreak\smallskip
  \noindent\kern12pt{\demofont@ #1.}\hskip8pt plus6pt minus4pt\ignorespaces
  }

\newif\ifextranotes	
\long\def\extranotes#1\endextranotes{\ifextranotes 
  \begingroup\par
    \leftskip=\parindent	\noindent
    {\tt Notes: }\ignorespaces#1\par\noindent
    {\tt\ End of Notes.}%
  \par\endgroup
\fi}

\define\http/{http:/{\kern-1pt}/}

\define\footnotesymlist{\yen\dag\ddag} 

\define\rotatesym#1#2.{%
  \global\let\holdsym=#1%
  \global\def\footnotesymlist{#2#1}
  }

\define\jkfootnote{%
  \expandafter\rotatesym\footnotesymlist.%
  \footnote"$^\holdsym$"%
  }

%

%


%

%
\define\ir#1{#1^\circ} 
%
\define\bd#1{{#1}^\partial}

\define\fnc#1,#2,#3.{{#1\:\mkern-1.5mu#2\mkern-2.5mu\to\mkern-2.5mu#3}}
\define\fnci#1,#2,#3.{{#1\:\mkern-1.5mu#2\mkern-2.5mu\hookrightarrow\mkern-2.5mu#3}}
\define\fncs#1,#2,#3.{{#1\:\mkern-1.5mu#2\mkern-2.5mu\twoheadrightarrow\mkern-2.5mu#3}}
\define\fto#1,#2.{{#1\mkern-2.5mu\to\mkern-2.5mu#2}} 
\define\ftoi#1,#2.{{#1\mkern-2.5mu\hookrightarrow\mkern-2.5mu#2}} 
\define\ftos#1,#2.{{#1\mkern-2.5mu\twoheadrightarrow\mkern-2.5mu#2}} 
%
\define\toitself{\looparrowleft} 
\define\fnctoitself#1,#2.{{#1\:\mkern-1.5mu#2\mkern-2.5mu\toitself}}


\define\setn#1{\lbrace#1\rbrace} 
\define\setg#1{\bigl\lbrace#1\bigr\rbrace}
\define\setG#1{\Bigl\lbrace#1\Bigr\rbrace}
\define\setd#1{\biggl\lbrace#1\biggr\rbrace} 
\define\setA#1{\left\lbrace#1\right\rbrace}

\define\setdef#1#2{\lbrace#1\mid#2\rbrace}
\define\setdefg#1#2{\bigl\lbrace#1\bigm|#2\bigr\rbrace}
\define\setdefG#1#2{\Bigl\lbrace#1\Bigm|#2\Bigr\rbrace}
\define\setdefd#1#2{\biggl\lbrace#1\biggm|#2\biggr\rbrace}
\define\setdefA#1#2{\left\lbrace#1\bigm|#2\right\rbrace}

		\message{   . . . leaving macros.tex }


\catcode`\@=\active

\input tiling.macros.ams


\def \eepacked {1}
\def \eetiled {2}
\def \eealgorithm {3}
\def \eetilepack {4}
\def \eethm {5}
\def \eeext {6}
\def \eec {7}
\def \eetable {8}
\def \eephenomenon {9}
\def \eeguess {10}
\def \eedede {11}
\def \eegoodexpr {12}
\def \eecoef {13}
\def \eearchetype {14}

\long\def\jept#1{{\rm({\ept #1})}}


\define\js#1{\lbrace#1\rbrace}
\define\jsd#1{\biggl\lbrace#1\biggr\rbrace}
\define\jsg#1{\bigl\lbrace#1\bigr\rbrace}
\define\jsG#1{\Bigl\lbrace#1\Bigr\rbrace}
\define\jsA#1{\left\lbrace#1\right\rbrace}
  \define\jS#1{\lbrace\!\lbrace#1\rbrace\!\rbrace}

\define\jp#1{(#1)}
\define\jpd#1{\biggl(#1\biggr)}
\define\jpg#1{\bigl(#1\bigr)}
\define\jpG#1{\Bigl(#1\Bigr)}
\define\jpA#1{\left(#1\right)}
  \define\jP#1{(\!(#1)\!)}

\define\jb#1{[#1]}
\define\jbd#1{\biggl[#1\biggr]}
\define\jbg#1{\bigl[#1\bigr]}
\define\jbG#1{\Bigl[#1\Bigr]}
\define\jbA#1{\left[#1\right]}
  \define\jB#1{[\![#1]\!]}

\define\ja#1{\langle#1\rangle}
\define\jag#1{\bigl\langle#1\bigr\rangle}
\define\jaG#1{\Bigl\langle#1\Bigr\rangle}
\define\jaA#1{\left\langle#1\right\rangle}
  \define\jA#1{\la\!\la#1\ra\!)} 

\define\jf#1{\lfloor#1\rfloor}
\define\jfg#1{\bigl\lfloor#1\bigr\rfloor}
\define\jfG#1{\Bigl\lfloor#1\Bigr\rfloor}
\define\jfA#1{\left\lfloor#1\right\rfloor}

\define\jc#1{\lceil#1\rceil}
\define\jcg#1{\bigl\lceil#1\bigr\rceil}
\define\jcG#1{\Bigl\lceil#1\Bigr\rceil}
\define\jcA#1{\left\lceil#1\right\rceil}

\define\jn#1{\|#1\|}
\define\jnd#1{\biggl\|#1\biggr\|}
\define\jng#1{\bigl\|#1\bigr\|}
\define\jnG#1{\Bigl\|#1\Bigr\|}
\define\jnA#1{\left\|#1\right\|}

\define\jv#1{|#1|}
\define\jvd#1{\biggl|#1\biggr|}
\define\jvg#1{\bigl|#1\bigr|}
\define\jvG#1{\Bigl|#1\Bigr|}
\define\jvA#1{\left|#1\right|}

\define\jcv#1{|#1|}
\define\jcvd#1{\biggl|#1\biggr|}
\define\jcvg#1{\bigl|#1\bigr|}
\define\jcvG#1{\Bigl|#1\Bigr|}
\define\jcvA#1{\left|#1\right|}

\define\card{\#} 

\define\jfo#1{ --\ignorespaces#1\ignorespaces-- } 
\define\jco#1{\unskip, \ignorespaces#1\ignorespaces,} 
\define\jq#1{``#1''}	\define\jQ#1{`#1'} 


  
\def\titlename{A change-of-coordinates from Geometry to Algebra,
		\\ applied to Brick Tilings}
\hcorrection{.25truein}	\pagewidth{6truein}
\vcorrection{0in}  \pageheight{8.5truein}
\baselineskip=14pt
\TagsOnLeft
  \nofootline
 \NoBlackBoxes
\extranotesfalse

\let\put=\relax \let\pf=\relax
\define\dpotwo/{diminished powers-of-two}
\define\logand/{{A{\kern-1.8pt}N{\kern-1.8pt}D}}
\define\logor/{{O{\kern-1.8pt}R}}

\define\FIR{1} \define\SEC{2} \define\THI{3} \define\LAS{\nproto}

\define\tables#1 {\message{OUTPUTTING Table}\nobreak\vskip\smallskipamount
  \bgroup\ept\narrower\narrower\narrower
  \noindent\hskip-12pt{\sc Tables #1\unskip\ \ }}
\def\endtables{\endgraf\penalty30
  \vskip\smallskipamount\egroup
  \noindent
  }

\redefine\footnotesymlist{\yen\dag\ddag}
\redefine\seq#1#2{{#1}^{\ja{#2}}}
\let\card=\jv
\let\tdm=\tau 
\define\nbrka{{\operatorname{\#Bricks}^\brka}}
%

\def\huchapternumber{0}	
\long\def\huchapter#1#2{\messageln{CHAPTER #1 #2}\goodbreak\bigskip 
  \vskip0pt plus 30pt \goodbreak
  \xdef\huchapternumber{#1}%
  \xdef\huchaptertext{#2}
  \centerline{\twelvepoint\sc \S #1 \ \ \ignorespaces #2}%
  \nobreak\medskip\noindent\tenpoint\ignorespaces
  }


\define\brkw{{\ss W}}
\define\brkx{{\ss X}}
\define\brky{{\ss Y}}
\define\brkz{{\ss Z}}
%

\define\latdede{{\Cal L}} 
\define\latprod{\Lambda}  
\define\stdprodlat{{\latprod\jb{\nproto,\dm}}} 
\define\iilat{{\latprod\jb{\ii,\ii}}} 

\define\ltw{\text{\tt w}} 
\define\ltx{\text{\tt x}}
\define\lty{\text{\tt y}}
\define\ltz{\text{\tt z}}

\define\alf{{\operatorname{Alf}}}	
\define\env#1{e^{#1}}		



\define\kbool{{Kin1}} 
\define\kk{Ke\&Ki} 


%

%
\topmatter
\hbox{}\endgraf\vskip -1\bigskipamount
\title \titlename \endtitle

\author\ept
\hbox{}\endgraf\vskip -2\bigskipamount
  Jonathan L.~King
\endauthor

\affil\ept
  \emf{University of Florida, Gainesville 32611-2082, USA} \;
\; {\tt squash\@math.ufl.edu}%
\; {\ss Webpage	http://www.math.ufl.edu/$\sim$squash/}%
\endaffil



\abstract
  A proof is sketched of the Polynomial Conjecture of the author
\jp{circulated as preprint \cite{\kbool},
\jq{Brick Tiling and Monotone Boolean Functions}%
},
which says
that the family of minimal tilable-boxes grows polynomially with dimension.
An important ingredient of the argument is translating the problem from its
fi\-nite-\dimal/ geometric framework to the algebraic setting of an
in-fi\-nite-\dimal/ lattice.
\endabstract
\endtopmatter

\document
\hbox{}\endgraf\vskip -3.6\bigskipamount
\huchapter1{Ingress}
What connection could there possibly be between packing boxes by $\nproto$
shapes of bricks, and the number of \logand//\logor/ logic circuits having $\nproto$
Boolean inputs?
\par
  Several years ago I found an algorithmic solution to a tiling problem,
aspects of which
\jco{it turned out}
had been solved more than two decades earlier.
\jkfootnote{This is not uncommon in tiling theory, whose literature-of-record
runs the gamut from technical research journals to puzzle books.}
  A sudden gust of serendipity, in the guise of Neil Sloane's {\tt superseeker}
program, led from my quantitative results to numerical evidence for a
\jq{Polynomial Conjecture} on the growth of
complexity \jp{rank} of the tiling space as a function of dimension.  In turn,
the conjecture
\jfo{henceforth abbreviated \jq{PC}}
led  to algebraic questions involving the \de/ sequence of integers
\jp{defined in \S2, along with PC}.
During a sabbatical year at the University of Toronto, a fruitful
collaboration with computer scientist Hugh Redelmeier led to additional numerical
support for the conjecture, then to a computer-assisted proof for $\nproto=5$.
In turn, this gave insight into the algebraic structure of the problem,
eventually culminating in a \jp{computerless} proof of the \pc/.

\specialhead
  Packings, Tilings  \& Algebra
\endspecialhead
  Packing-type problems are arguably among the most ancient of combinatorial
conundrums.
In recent times, several types of overtly algebraic methods have been used to study
packings/tilings.
\roster
\item"$\bullet$"
  Group theory, in the form of symmetry groups of tessellations and of
crystals.
\item"$\bullet$"
  Combinatorial group theory, e.g, \cite{Co\&La} \cite{Thurst} \cite{Propp}.
\item"$\bullet$"
  Commutative algebra, e.g, \cite{Bar{1,2}}.
\goodbreak
\item"$\bullet$"
  \jq{Dehn invariant}  and related Tensor Algebra methods, e.g,
\cite{Dehn} \cite{Lac\&Sze} \cite{Fr\&Ri} \cite{Gal\&G} \cite{\kk1,2}.
\endroster
\goodbreak\noindent
The tiling problem of the present paper is the first, to my knowledge, which
seems to require a smidgeon of Lattice Theory.
  In this note, I will sketch the passage from the Geometry to the Algebra
\jp{from tilings to lattices}
and show how PC reduces to a \jq{finiteness certificate} which can be verified
by computer.
\par
  The current {\S1} defines brick tilings and the \jq{rank} of a set of
protobricks.
{\S2}~states PC, illustrates how rank can be computed, and gives a brief
introduction to distributive lattices and the \de/ sequence.
In {\S3}, PC is restated in the lattice setting, and the finiteness certificate
is described along with a sketch of  ideas employed in the proof.
The technical lattice-algebraic demonstration of PC will appear in \cite{Kin5}.
Geometric information on brick tilings/packings appears in
preprints \cite{\kbool}~and~\cite{Kin2}.  Lastly, {\S4} lists open questions.

\subhead
  Brick-Packings/Tilings
\endsubhead
  Published in innumerably many puzzle books is this chestnut:
\emfs{Can the $8\times8$ chessboard, minus its North-East and South-West
corner squares, be packed by \jp{thirty-one} dominos?}
\jp{The dominos can be placed in both the $1\ct2$~and $2\ct1$ orientations.}
\par
  We are \jq{born knowing} that there is no such packing:  The two removed squares
have the same color, yet each orientation of a domino must cover both a black and a
white square.
  Even should we allow ourselves to place positive \emf{and negative}
copies of dominos
\jp{{\ept defined at~(\eetiled)}}, 
this \jq{coloring argument} still precludes a tiling.
Depending on the shapes of the \jq{proto-tiles}, coloring ideas sometimes give
an IFF-condition for whether a specified region is tilable.

\subsubhead
  Bricks
\endsubsubhead
  A $\dm$-dimensional \df{brick} $\brkb$ is a $\dm$-tuple
	$$
\brkb \;=\;
\slb_1 \times \cdots \times \slb_\dira \times \cdots \times \slb_\dm
\,,	 $$%
where each \df{sidelength}~$\slb_\dira$ is a positive integer.
\jkfootnote{Brick-tiling questions which permit non-integral sidelengths are
discussed in \cite{Lac\&Sze} \cite{Fr\&Ri} \cite{\kk{1,2}}.}
  We will identify each
{$\dm$-brick} $\brkb$ with a product of half-open intervals,
	$$
\brkb \;=\;
{\rco0,\slb_1.}\times \cdots\times{\rco0,\slb_\dm.}
\; ,	 $$%
a subset of Euclidean space~$\reald$.
Translating our brick by a vector, $\brkb+{\vec w}$, gives the set of all sums
${\vec y}+{\vec w}$ for ${\vec y}\in\brkb$.
Agree to use $\brka,\brkb,\brkc,\brkt$ to name bricks.  A lowercase
letter denotes the corresponding sidelengths, e.g,
	$$
\brka=\sla_1 \times \cdots \times \sla_\dm
\quad\text{and}\quad
\brkt=\slt_1 \times \cdots \times \slt_\dm
\,.	 $$%
A \df{box} is another name for a brick; the latter are used to pack/tile the
former.
I am interested in when a specified box~$\brkt$ \jp{the target} can be packed
or tiled by translates of copies of bricks in a specified finite set
	$$
\proto \;=\; \jsg{{\seq\brkb\FIR},{\seq\brkb\SEC},\dots,{\seq\brkb\LAS}}
	 $$%
called the set of \dfw{protobricks}.

\subsubhead
  Definitions: Packing \& Tiling
\endsubsubhead
  For a subset $S\subset\reald$, the indicator function $\ind{S}$ is~$1$ at
those points ${\vec y}$ in~$S$, and $\ind{S}\jp{\vec y}$ is $0$ on the
complement~${\reald\setdiff{S}}$.
In order to show the connection between the problems considered in
\cite{Bar{1,2}} and \cite{\kk1,2}, I define \jq{tilable} a touch more
generally than is needed in the present paper.
\par
  Given a set $\proto$ of protobricks in~$\reald$, a
box $\brkt$ is \df{packable} if
	$$
\ind\brkt \;=\; \sum_{\brkh\in\bsetb}\ind{\brkh} \, ,
\tag\put\eepacked
	 $$%
for some finite collection $\bsetb$ of translates of the protobricks.
{Figure\,$\eepacked'$} exhibits a packing, by $\proto=\js{\brka,\brkb,\brkc}$.
\par
\lfigure
  {\epsfxsize=2.7in\epsfbox{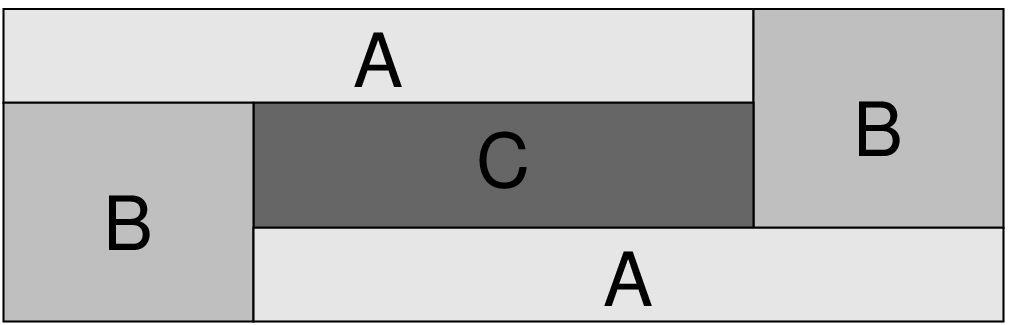}}
  {1.5in}
  {$\eepacked'$}
  Two copies of rectangle~$\brka$, two of $\brkb$ and a single $\brkc$ pack
the $34\ct11$ rectangle, where
the protobricks are
	$$
\spreadlines{-1\jot}
    \matrix
\brka\:\quad & 25 &\times &3
    \\
\brkb\:\quad & 9  &\times &8
    \\
\brkc\:\quad & 16 &\times &5
    \endmatrix
	 $$%
\endlfigure
\par
  For tiling, I allow weights from an arbitrary commutative monoid
$(\Gamma,{+},0)$, with a distinguished non-zero element $1\in\Gamma$.  Say
that a box~$\brkt$ is \df{$\Gamma$-tilable} \jco{by~$\proto$} if
there exists a finite collection $\bsetb$ of protobrick translates as well as
coefficients 
\jkfootnote{%
  For brick tiling, both \cite{Bar2, thm\,2.1} and \cite{\kbool,
Equality\,Thm} show that $\complex$-tilability is equivalent to
$\integers$-tilability.  In contrast, \cite{Bar2, P.14} has an example of a
box which can be $\rationals$-tiled by certain \emf{polyominos}, but cannot be
$\integers$-tiled by them.
} 
$\gamma_\brkh\in\Gamma$, for $\brkh$ in~$\bsetb$, such that
	$$
\ind\brkt \;=\; \sum_{\brkh\in\bsetb} \gamma_\brkh\ind{\brkh} \, .
    \rlap
	{$\qquad\qquad\foldedtext\foldedwidth{0.8in}{\ept
          \llap{(}The addition takes place in $\Gamma$.\rlap{)}%
	  }$
	}
\tag\pf\eetiled
	 $$%
When $\Gamma$ is {the additive group of integers},
say simply that $\brkt$ is \df{tilable} by~$\proto$.
As an
illustration, consider the protobrick set consisting of three rectangles
$\brka=3\ct8$, $\brkb=4\ct5$ and $\brkc=7\ct3$.
Our target is~$\brkt\idtb3\ct1$.
\rfigure
  {\epsfysize=1.4in\epsfbox{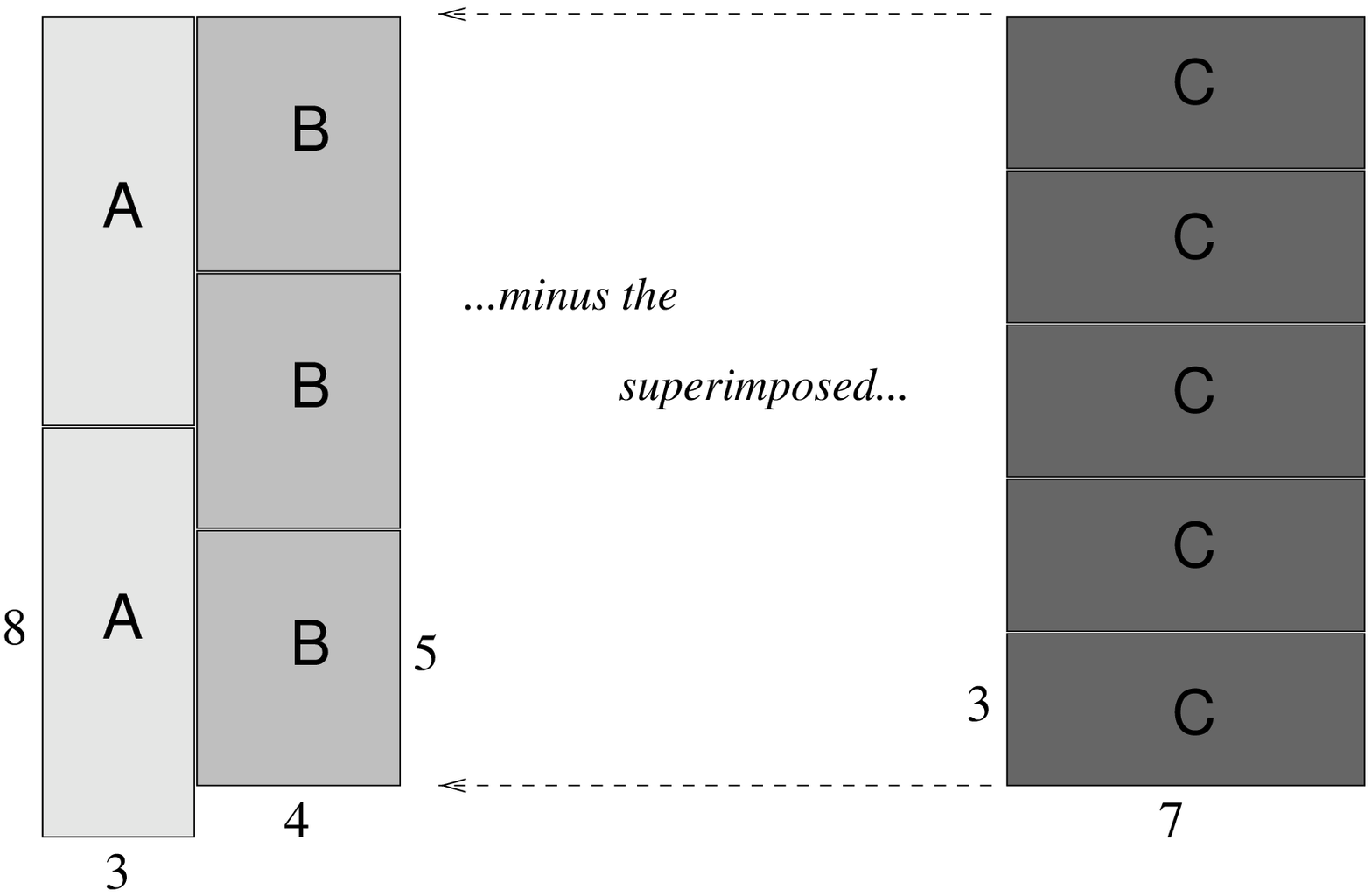}}
  {1.85in}
  {$\eetiled'$}
  Although it is drawn on the righthand side of the image, the
\jq{$\brkc$\,clump} is to be superimposed on the \jq{$\brka\brkb$\,clump}, with
left (and hence right) edges aligned.  Note that {subtracting} the
\jq{$\brkc$\,clump} from the \jq{$\brka\brkb$ clump} leaves
a tiny $3\ct1$~rectangle at the bottom of the lowermost~\jq{$\brka$}.
\endgraf
    Thus, by using $2$\,copies of~$\brka$ and $3$\,copies of~$\brkb$
and $(-5)$\,copies of~$\brkc$, we tile~$\brkt$.
\endrfigure
  This figure depicts a way to tile~$3\ct1$ by the proto-set~$\setn{\brka,\brkb,\brkc}$.
Certainly $\brkt$ cannot \jco{however} be \emf{packed} by these protobricks.

\subsubhead\nofrills
  Why Tiling?\usualspace
\endsubsubhead
  Trivially, there is an algorithm which is exhaustive
\jfo{indeed, exhaust\-ing} for determining whether a given $\brkt$ is packable
\jept{by a fixed protoset~$\proto$}.  At first glance, one might think that
tilability of~$\brkt$ is more difficult to ascertain since, potentially,
there are infinitely many collections $\bsetb$, in~($\eetiled$), to consider.
It transpires that the opposite is true.  There is an analogy
	$$
    \align
\text{Packings} &\hookrightarrow \text{Tilings}
    \\
\text{Semigroup} &\hookrightarrow \text{Group}
    \endalign
	 $$%
between studying a semigroup by embedding it in a group, and studying packings
by first understanding the \jp{larger} space of tilings.  And groups are
easier to understand than
semigroups, as everybody knows\dots

\subhead
  The partial order on Bricks
\endsubhead
  The goal of \cite{\kbool} was to find a fast algorithm, with $\proto$ fixed, for
determining whether a target box~$\brkt$ is packable or tilable.
\par
  On the space of $\dm$-bricks there is a natural partial order \jq{$\brkle$} of
packability.
For integers, use $\sla\slle\slb$ for \jq{\emfs{$\sla$~divides $\slb$}} and
$\slb\slge\sla$ for \jq{\emfs{$\slb$~is a multiple of~$\sla$}}.
Say that $\brkb$ \df{\parapack/s~$\brkt$} \jept{or \df{divides~$\brkt$}}
\jco{written $\brkb\brkle\brkt$}
if translates of~$\brkb$ can pack~$\brkt$.  That is, {$\brkb\brkle\brkt$} iff
	$$
\text{For each $\dira=1,\dots,\dm$ direction: \;
   sidelength $\slb_\dira$ divides $\slt_\dira$.%
}
	 $$%
Now consider a collection $\bsetb$ of boxes which is an \dfw{up-set} in the
partial order,
	$$
\display
  Each box which is a multiple of some $\bsetb$-brick, is necessarily itself a $\bsetb$-brick.
\enddisplay
	 $$%
Each up-set $\bsetb$ is determined by its family of minimal elements
\jept{w.r.t\. the~$\brkle$ order}.  Writing this family
as $\minmal\jp\bsetb$, we
have that ${\brkt\in\bsetb}$ iff
	$$
\exists{\brkb\in{\minmal\jp\bsetb}}
\,\;\text{with}\;\,
{\brkb\brkle\brkt} \,.
\tag\put\eealgorithm
	 $$%
Thus both $\pkable\jp\proto$
\jfo{the set of $\proto$-packable boxes}
and $\tlable\jp\proto$
\jfo{the tilable boxes}
are determined by their sets of minimal elements, respectively.
\par
  My purpose in \cite{\kbool}
\jco{partially successful}
was to find finite descriptions of $\pkable\jp\proto$ and $\tlable\jp\proto$
which allowed an {efficient} test for membership.
Alas for packing, the minimal set doesn't work; typically
$\minmal\jpg{\pkable\jp\proto}$ is infinite.
%
%
\par
  Happily, several authors
%
\jkfootnote{%
  Although they do not discuss the minimal set of~$\tlable\jp\proto$, {Katona
\& Sz\'asz (1971)} give a criterion for a box being $\proto$-tilable
\jp{and $\proto$-packable, once the sidelengths are large enough}
under the
assumption that $\proto$ comprises all $\dm!$~orientations of bricks in a
brick-set.
\endgraf
  Barnes, in two seminal papers (1982), uses ideals over polynomial rings to
develop a general algebraic criterion for a polyomino to be tilable by other
polyominos.  Results {(2.1)} and {(2.4)} of \cite{Bar{2}} imply our
{theorem~\eetilepack}, here.
\endgraf 
  These three papers do not address computability.}
proved versions of the important result that the set
${\bsetm\jp\proto}\idtb{\minmal\jpg{\tlable\jp\proto}}$
is finite, and that sufficiently large boxes are $\proto$-packable if and only
if they are $\proto$-tilable.
\jb{{\ept See the end of~{\S2} for an example computation of~$\bsetm\jp\proto$.}}

\proclaim{Theorem \pf\eetilepack}
 $\bsetm\jp\proto$ is finite and is computable.  Furthermore, there is a computable
integer $\bndpackable=\bndpackable\jp\proto$ so that, whenever $\brkt$ is a box
whose sidelengths each exceed $\bndpackable$, then:
\quad
$\text{$\brkt$ tilable}\implies\text{$\brkt$ packable.}$
\endproclaim

This theorem yields an algorithm for testing whether a candidate box
$\brkt$ is tilable:
\;\emfs{Does $\brkb\brkle\brkt$, for some~$\brkb$ in ${\bsetm\jp\proto}$?}
\quad
  Letting $n$ denote the number of bits needed to describe $\brkt$, this
algorithm runs in linear time~$O(n)$.

\specialhead
  Computing rank
\endspecialhead
  Is the algorithm practical?  Well\dots, this all depends on the magnitude of
the constant in the $O(n)$~algorithm.
\cite{\kbool} called the cardinality of~${\bsetm\jp\proto}$
the \df{rank} of~$\proto$, showed it bounded by a pure function of
$\nproto$~and $\dm$
\jept{the number of protobricks and their dimension},
and produced two algorithms for computing it.  Here is the first algorithm:
\par
  Given bricks $\brka,\brkb,\dots,\brkc$, we can use them to tile a box
$\brkt$ built as follows.
Let
	$$
    \align
g_1 \;&\idtb\; \gcd\js{\sla_1,\slb_1,\dots,\slc_1},
\quad\text{and for each other direction:}
    \\
\ell_\dirb \;&\idtb\; \lcm\js{\sla_\dirb,\slb_\dirb,\dots,\slc_\dirb},
\quad\text{for $\dirb=2,3,\dots,\dm$.}
    \endalign
	 $$%
Then
{\,$\brkt\idtb{g_1\times\ell_2\times\cdots\times\ell_\dm}$\,}
is tilable by collection $\js{\brka,\brkb,\dots,\brkc}$.  To see this, note
that $\brka$ \parapack/s the \jq{slab}
	$$
\brka'\;\idtb\;{\sla_1\times\ell_2\times\cdots\times\ell_\dm}
\,.	 $$%
And $\brkt$ is tiled by slabs $\js{\brka',\brkb',\dots,\brkc'}$ in the same way that
integer~$g_1$ is an integral linear-combin\-a\-tion of integers
$\js{\sla_1,\slb_1,\dots,\slc_1}$.
We call this the \df{combine} operation, and write
$\brkt=\comb_1\jpg{\js{\brka,\brkb,\dots,\brkc}}$, the combine in direction~$1$.
\par
  More generally, given a brick-set $\bseta$ and direction~$\dira$, let
$\prjt\bseta\dira$ denote the \emf{set} $\setdef{\slb_\dira}{\brkb\in\bseta}$
of $\dira$th sidelengths.  Then $\comb_\dira\jp\bseta$ is the brick
$\slt_1 \ct \cdots \ct \slt_\dm$, where
	$$
    \align
\slt_{\dira} &\idtb\gcd\jpg{\prjt\bseta\dira}
\quad\text{and, for each direction $\dirb\neq\dira$:}
    \\
\slt_{\dirb} &\idtb\lcm\jpg{\prjt\bseta\dirb} \,.
    \endalign
	 $$%
\par
  It turns out that iterating all possible Combines is powerful enough to
generate all of $\bsetm\jp\proto$.
Define the \dfw{$\dira$th extension of~$\proto$} to be the set of bricks
	$$
\ext_\dira(\proto) \;\idtb\;
\setdefg{\comb_\dira(\bseta)}{\text{$\bseta$~is a non-void subset of~$\proto$}} \,.
	 $$%
It is not difficult to see that each two of the $\ext$ operators commute,
and each is idempotent. So~(\eeext\?a) \jco{below} is the set of all boxes
that can be made by means of the $\comb$ine operation.
Moreover,

\proclaim{Theorem~\pf\eethm  {\;\ept\cite{\kbool, Equality Thm}}}
  The brick-set $\bsetm\jp\proto$ equals the set of minimal bricks
{\ept\sl\jp{w.r.t\.~divisibility}}
of
	$$
\ext_\dm\jpG{\ext_{\dm-1}\jpg{\dots\ext_2\jpg{\ext_1\jp{\proto}}\dots}} \,.
\tag\pf\eeext\?a
	 $$%
\endproclaim

As a corollary, we get this daunting bound on the rank of~$\proto$.
	$$
\rank\jp\proto
\;\le\;
2^{2^{2^{\vdots^{2^\nproto}}}
    \rlap{\quad\vbox{\hsize=1.3in\noindent\ept
          \llap{(}a tower of~$\dm$\,many \endgraf\noindent exponentiations\rlap{)}} }
  }
\tag\eeext\?b
	 $$%
\par
  On the one hand \jfo{perhaps unexpectedly} formula~(\eeext\?a) gives a
workable algorithm for computing the set ${\bsetm\jp\proto}$ of minimal
tilable-boxes.  For as the algorithm progressively generates bricks, we can
discard bricks when they become divisible by a later-generated brick.
\par
  On the other hand, bound~(\eeext\?b) is laughably too large.  Other
arguments in \cite{\kbool} give a smaller bound of
$\rank\jp\proto\le{\dm^{\dm^\nproto}}$.  The corresponding algorithm, however,
which this smaller bound engenders, typically runs more slowly than that
from~(\eeext).

\subhead
  Two examples
\endsubhead
  The rank of~$\proto$ can be smaller or larger than cardinality~$\card\proto$.
\par
  Proto-set $\proto=\js{\brka,\brkb,\brkc}$
of Figure~$\eepacked'$ tiles~$34\ct11$.  What is the rank
of~$\proto$?  Evidently $\comb_1$ produces these bricks,
\vskip -\medskipamount
	$$
    \align
\comb_1\js{\brka,\brkb} \;&=\; 1\times(3\cdot 8)
    \\
\comb_1\js{\brkb,\brkc} \;&=\; 1\times(8\cdot 5)
    \\
\comb_1\js{\brkc,\brka} \;&=\; 1\times(5\cdot 3)
\, .
   \endalign
	 $$%
Applying $\comb_2$ to this family of three bricks yields the $1\ct1$~brick.
Thus $\tlable\jp\proto$ is the set of \emf{all} boxes.  So $\rank\jp\proto$
is~$1$.

\smallskip
  As a second example,
let $\brka$ be $2\ct3\ct7$ and let $\proto\idtb\jsg{\brka,\brka',\brka''}$,
where each stroke means to rotate the sides by one position;
$\brka'=3\ct7\ct2$ and $\brka''=7\ct2\ct3$.
Necessarily, the set ${\minmal\jpg{\tlable\jp\proto}}$  \jco{the minimal tilable-boxes}
is rotation invariant.  It comprises these five bricks
	$$
  \matrix
\format \l   &
  \c &	\;\c\; &		\c  &	\;\c\; &	\c
  \\
\brka 				 & \,:\quad
  2 &	\times &	{3} &	\times &		7
  \\
\brkb		\idtb\comb_1\js{\brka,\brka'} & \,:\quad
  1 &	\times &	{(3\cdot 7)} &	\times & 	{(7\cdot2)}
  \\
\dot\brkb	\idtb\comb_1\js{\brka',\brka''} & \,:\quad
  1 &	\times &	{(7\cdot 2)} &	\times & 	{(2\cdot3)}
  \\
\ddot\brkb	\idtb\comb_1\js{\brka'',\brka} & \,:\quad
  1 &	\times &	{(2\cdot 3)} &	\times & 	{(3\cdot 7)}
  \\
\brkc	\idtb\comb_2\js{\brkb,\dot\brkb,\ddot\brkb} & \,:\quad
  1 &	\times &	{1} &	\times	& 	{(2\cdot 3\cdot 7)}
  \endmatrix
	 $$%
and their rotates.  Thus $\rank\jpg{\js{\brka,\brka',\brka''}}=15$.

\subhead
  Simplifying $\ext$~notation
\endsubhead
  Given a finite set $\zseta$ of directions \jp{positive integers}, let
$\ext_\zseta$ mean
	$$
\ext_{\dira_1} \circ \ext_{\dira_2} \circ \dots \circ \ext_{\dira_\tdm}
\,,	 $$%
where $\dira_1,\dots,\dira_\tdm$ is some enumeration of~$\zseta$; this is
well-defined since all the $\ext$ operators commute.   Henceforth, write
$\ext_\dm\jpg{\cdots{\ext_1\jp{\proto}}\cdots}$
as
$\ext_{\js{1,\dots,\dm}}\jp\proto$, or just as
$\extdfold\jp\proto$.  When $\zseta$ is empty,
$\ext_\emptyset\jp\proto$ is~$\proto$.

\huchapter2{Polynomial Conjecture}
  Since the tiling-rank of a set of $\nproto$ many $\dm$-bricks is bounded by
a function of $\nproto$~and $\dm$, and since rank is essentially the constant
in the linear-time algorithm tilability test, one naturally wishes to study
the \df{maxrank} function~$\maxrank$:
	$$
\display
  $\maxrank\jp{\nproto,\dm}$ is the maximum, as $\proto$ ranges over
all $\nproto$-sets of $\dm$-dimensional bricks, of $\rank(\proto)$.
\enddisplay
	 $$%
For each $\nproto,\dm$~pair, there is a straightforward method to construct a
\jq{worst case} brick-set $\proto$ whose rank is $\maxrank\jp{\nproto,\dm}$.
Letting
$\slb^{\ut\dm}$
denote the $\dm$-cube $\slb\times\slb\times\cdots\times\slb$, it turns out that a worst
case~$\proto$ can be built from \emf{cubes},
	$$
\proto\;=\;
\jsG{{\jp{\seq\slb\FIR}^{\ut\dm}}, \; {\jp{\seq\slb\SEC}^{\ut\dm}}, \;
  \dots, \; {\jp{\seq\slb\LAS}^{\ut\dm}} } \,.
\tag\put\eec\?a
	 $$%
Moreover, the collection of $\nproto$~sidelengths
$\seq\slb\FIR, \seq\slb\SEC,\dots,{\seq\slb\LAS}$,
can be chosen to depend only on~$\nproto$, and not on dimension.  For
$\nproto=3$, here is one such collection.
	$$
    \align
{\seq\slb\FIR} \;&\idtb\;
2^1 \cdot 3^1 \cdot 5^2 \cdot 7^2 \cdot 11^3 \cdot 13^3
    \\
{\seq\slb\SEC} \;&\idtb\;
2^2 \cdot 3^3 \cdot 5^1 \cdot 7^3 \cdot 11^1 \cdot 13^2
    \\
{\seq\slb3} \;&\idtb\;
2^3 \cdot 3^2 \cdot 5^3 \cdot 7^1 \cdot 11^2 \cdot 13^1
\,.
    \endalign
	 $$%
Reading the exponents down the columns, we see each of the
$6$~permutations of~$\js{1,2,3}$.
\par
  More generally, let
$\setdef{p_\perm}{\perm\in\text{Perms}}$
be the first $\nproto!$~prime numbers, indexed by the $\nproto!$ permutations
of~$\js{\FIR,\dots,\LAS}$.

\proclaim{Proposition~{\eec\?b}  {\;\ept\cite{\kbool, Max-rank Proposition}}}
  For $n=\FIR,\dots,\LAS$, let
	$$
\seq\slb{n} \;\idtb\; \prod_{\perm\in\text{Perms}} p_{\!\perm}^{\,\perm\jp{n}}
\,.	 $$%
Then proto-set~$(\eec\?a)$ has rank equal to the maxrank value $\maxrank\jp{\nproto,\dm}$.
\endproclaim

So maxrank values can now be computed.  A program I wrote in Common Lisp
calculated the table below.
%
\nobreak
\medskip
\setbox0\hbox{$
\matrix
  \format \l\quad&&\quad\r
\\
  \relax  & \llap{$\dm\!\rightarrow$\;}\dimnum2  & \dimnum3  & \dimnum4
  & \dimnum5  & \dimnum6  & \dimnum7  & \dimnum8
\\
  \bold1  &                     1 &   1 &    1 &    1 & 1   & 1   &  1
\\
  \bold2  &                     4 &   5 &    6 &    7 & 8   & 9   & 10
\\
  \bold3  &                    18 &  36 &   61 &   93 & 132 & 178 & 231
\\
  \bold4  &                   166 & 578 & 1372 & 2669 & 4590 & ?
\\
  \bold5  &		     7579 &   ?
\endmatrix
	 $}%
\centerline{%
  $\vcenter{\hbox{\llap{$\card\proto=\,$}\hbox{$\nproto$}}\hbox{$\;\downarrow$}}$
  $\vcenter{\box0}$
  }
\nobreak
\table{\put\eetable}
 {Maximal-rank values, $\maxrank(\nproto,\dm)$}.  The computer program omitted
\jp{serendipitously, as it later turned out}
the $\dm\!=\!\dimnum1$ column, which is trivially constant~$1$.
\endtable

\subhead
  Enter Computer Serendipity
\endsubhead
  Neil Sloane's venerable \jq{Handbook of Integer Sequences} was the standard
reference for looking up mystery sequences of integers.  As a wonderful
service to the mathematical community, he has now made available an electronic
version \jco{called {\tt superseeker}} which, upon receiving an email message
comprising some terms of a sequence, mails back a list of journal citations
where the sequence has been analysed.
\par
  The computed $\dm=\dimnum2$ column,
$1,4,18,166,7579$,
were the only numbers in the table not entirely mysterious to me.  The proof of~(\eec\?b)
showed that they were the \emf{\de/ numbers} \jept{defined below}, an explosively-growing sequence
of integers well-known to combinatorists.
In order to save a trip to the library, I used Sloane's program to get
citations for \de/'s sequence.
On a lark, I subsequently emailed
\jq{{\tt lookup 18 36 \dots 231}} \jco{the $\nproto=\bold3$ row}
to {\tt superseeker\@research.att.com}.  Knowing that the
\de/ sequence grew \emf{doubly}-exponentially with~$\nproto$, then, I was dumbfounded to receive
	$$
\vbox{\ept\parindent=4\parindent \parskip=0pt \tt
  From: superseq-reply\@research.att.com

  To: squash\@math.ufl.edu

  \smallskip
  Report on [ 18,36,61,93,132,178,231]:

  Many tests are carried out, but only potentially useful information

  (if any) is reported here.

  \smallskip
  TEST: IS THE k-TH TERM A POLYNOMIAL IN k?
  \par
  \hbox{}\hskip30pt
	  SUCCESS: k-th term is nontrivial polynomial in k of
  \par
  degree \quad   2

  \smallskip
  Polynomial is    18+29/2*k+7/2*k$\uparrow$2
  \par
}
	 $$%
\par
  Completely floored, I hastily emailed off what little I had for the
$\nproto=\bold4$~row, \,only to see
	$$
\vbox{\ept\parindent=4\parindent \parskip=0pt \tt
  Report on [ 166,578,1372,2669,4590]:
  \smallskip
  TEST: IS THE k-TH TERM A POLYNOMIAL IN k?

  \smallskip

  \hbox{}\hskip30pt
	  SUCCESS: k-th term is nontrivial polynomial in k of
  \par
  degree \quad   3

  \smallskip
  Polynomial is    166+784/3*k+261/2*k$\uparrow$2+121/6*k$\uparrow$3
}
	 $$%
This was certainly food for contemplation\dots --perhaps row~$\nproto$ was the
output of a degree $\nprotomo$ polynomial?!  Writing down the apparent
polynomials for $\nproto=\bold1,\bold2$ gave this list.
\define\bob{h}
	$$
    \aligned
\bob_{\bold1}(k)\;&\idtb\;
1
    \\
\bob_{\bold2}(k)\;&\idtb\;
4 + k
    \\
\bob_{\bold3}(k)\;&\idtb\;
18+ \tfrac{29}{2}k + \tfrac{7}{2}k^2 \,.
    \\
\bob_{\bold4}(k)\;&\idtb\;
166 + \tfrac{784}{3}k + \tfrac{261}{2}k^2 + \tfrac{121}{6}k^3 \,.
    \endaligned
	 $$%
\par
  The penny still had not dropped; I saw no pattern in this list.
Moreover, plugging $k=-1$ \jp{which corresponds to~$\dm=\dimnum1$} into
$\bob_{\bold3}$ did \emf{not} give the correct value of
$1=\maxrank\jp{{\bold3},{\dimnum1}}$, but rather gave~$7$.
\par
  But stay a moment -evaluating \emf{all} the polynomials
$\bob_{\bold1},\bob_{\bold2},\bob_{\bold3},\bob_{\bold4}$ at~$k=-1$, gave
$1,3,7,15$ --the \dpotwo/?!  And plugging in $k=-2$ yielded
$1,2,3,4$.  Hmm\!\dots
\par
At last the penny dropped.
	$$
\display
  Some phenomenon in {Table~$\eetable$} only kicked in at $\dm\ge{\dimnum2}$.
However, the phenomenon was naturally indexed from~$\dm={\dimnum0}$.
\enddisplay
\tag\put\eephenomenon
	 $$%
Shifting the polynomials back by~$2$
\jco{by letting $\poly_\nproto\jp\dm\idtb{\bob_\nproto\jp{\dm-2}}$}
gave the first $4$~lines of these next two tables.  \jept{The $5$th~lines were
computed later, by Hugh Redelmeier.}
	$$
    \aligned
\poly_{\bold1}(\dm)\;&\idtb\;
1
    \\
\poly_{\bold2}(\dm)\;&\idtb\;
2+\dm
    \\
\poly_{\bold3}(\dm)\;&\idtb\;
3 + \tfrac{1}{2!} \bigl[ \dm + {7}\dm^2 \bigr]
    \\
\poly_{\bold4}(\dm)\;&\idtb\;
4 + \tfrac{1}{3!}\bigl[ {-112}\dm + {57}\dm^2 + {121}\dm^3 \bigr]
    \\
\poly_{\bold5}(\dm)\;&\idtb\;
5 + \tfrac{1}{4!}\bigl[
  29898\dm - 81241\dm^2 + 48066\dm^3 + 3901\dm^4
\bigr]
\qquad{\text{\jp{Redelmeier}}}
    \endaligned
	 $$%
%
\nobreak
\medskip
\setbox0\hbox{\eightpoint$
\matrix
  \format \l\,\quad&&\quad\r
\\
  \relax  & \llap{$\dm\!\rightarrow$\;} \dimnum0 & \dimnum1 & \dimnum2  & \dimnum3  & \dimnum4
  & \dimnum5  & \dimnum6  & \dimnum7  & \dimnum8  & \dimnum9 & \dimnum{10} & \dimnum{11}
\\
  \bold1  &	1 &  1 &	1 &   1 &    1 &    1 &  1 & 1 &  1 & 1 & 1 & 1
\\
  \bold2  &	2 &  3 &	4 &   5 &    6 &    7 & 8 & 9 & 10 & 11 & 12 & 13
\\
  \bold3  &	3 &  7 &	18 &  36 &   61 &   93 & 132 & 178 & 231 & 291 & 358 & 432
\\
  \bold4  &	4 & 15 &	166 & 578 & 1372 & 2669 & 4590 & 7256 & 10788 & 15307 & 20934 & 27790
\\
  \bold5  &	5 & 31 & 7579 & 40517 & 120614 & 273540
  &&& {\llap{\jp{Redelmeier}}}
\endmatrix
	 $}%
\hbox{%
  $\vcenter{\hbox{\eightpoint$\nproto$}\hbox{$\,\downarrow$}}$
  $\vcenter{\box0}$
  }
\nobreak
\tables{\put\eeguess\?a \& \eeguess\?b}
   Polynomials $\poly_\nproto$ and their values $\poly_\nproto(\dm)$,
for~$\nproto\le\bold5$.
The zero-th column shows the naturals, the first exhibits \dpotwo/, and
column~$\dimnum2$ has \de/ numbers.
\endtables
With this much mathematical smoke in evidence, it was irresistable to
conjecture that there was some mathematical fire underlying it.

\proclaim{Polynomial Conjecture}
  As~$\dm$ ranges over $\ico2,\ii.$, the mapping $\dm\mapsto\maxrank(\nproto,\dm)$
is a polynomial
\jkfootnote{Every degree-$\jp\nprotomo$ polynomial which takes on integer values at integers
necessarily has coefficients of the form 
$q\!\bigm/\!{\jp\nprotomo!}$ where $q$ is integral.
This form of the coefficients will arise naturally
\jco{in~(\eecoef)}
from the proof of PC in~{\S3}.
} 
%
of degree~$\nprotomo$.
\endproclaim

  In order to get a handle on this conjecture, we study the
algebraic structure underlying~(\eeext\?a), which is that of a lattice.

\specialhead
  Lattices {\,\&\,} CIA {\!}operators
\endspecialhead
  A \df{lattice} is a poset ${\jp{\latdede, \le}}$ such that each pair
$\sla,\slb\in\latdede$ has a \emf{greatest lower bound} \jco{written
$\sla\wedge\slb$}
and a \emf{least upper bound}, $\sla\vee\slb$.  Letting
$\slc\idtb{\sla\wedge\slb}$, then,
$\slc\le\sla$ \& $\slc\le\slb$ and, if $\slc'$ is any other such element, then
$\slc\ge{\slc'}$.
\par
  Automatically, $\wedge$~is a \dfw{CIA operator}
\jfo{Commutative, Idempotent \jept{$\sla\wedge\sla=\sla$}, Associative}
and so is~$\vee$.
  Moreover, the lattice operations fulfill the absorption laws
	$$
\sla\vee{\jp{\sla\wedge\slb}} = \sla
\qquad\&\qquad
\sla\wedge{\jp{\sla\vee\slb}} = \sla \,.
	 $$%
\par
  Now consider the brick-set
$\bseta\idtb\extdfold\jp\proto$, for a fixed~$\dm\ge2$.
Observe that the collection $\prjt\bseta1$ of first sidelengths is a
lattice with respect to divisibility; here \jq{$\wedge$} is~$\gcd$ and
\jq{$\vee$} is~$\lcm$.  Indeed, in each direction~$\dira$,
	$$
\display
  $\jpG{\, {\prjt\bseta\dira} \,,\; \slle\;}$ is a distributive lattice,
\enddisplay
	 $$%
since
$\sla\vee\jp{\slb\wedge\slc}=\jp{\sla\vee\slb}\wedge\jp{\sla\vee\slc}$,
and
$\sla\wedge\jp{\slb\vee\slc}$ equals
$\jp{\sla\wedge\slb}\vee\jp{\sla\wedge\slc}$.
\par
  Each lattice $\prjt\bseta\dira$ is generated by the
\jp{at most} $\nproto$~numbers $\prjt\proto\dira$.  Thus each of these
lattices is a homomorphic image of $\latdede\jb\nproto$
\jco{the \df{\de/ lattice}} which is the \emf{free} distributive lattice on
$\nproto$~generators.

\subhead
  Picturing $\latdede\jb\nproto$
\endsubhead
  It is convenient to describe $\latdede=\latdede\jb\nproto$ as the lattice of
non-decreasing Boolean functions of $\nproto$ Boolean variables.  To this end,
we will write \jq{$\vee$} as addition
\jco{logical \logor/}
and write \jq{$\wedge$} as
multiplication, {logical~\logand/}.
\rfigure
  {\epsfysize=3in\epsfbox{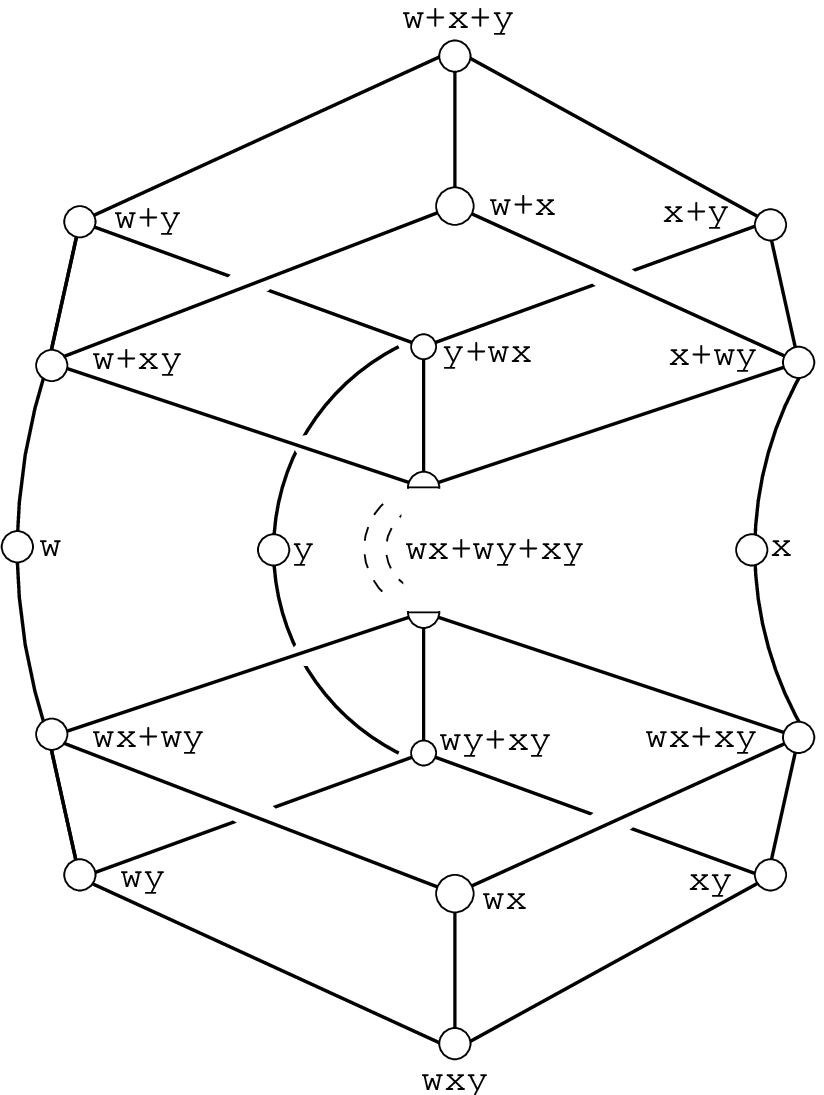}}
  {2.10in}
  {\put\eedede}
 The \de/ lattice $\latdede\jb3$, here generated by the three symbol
alphabet $\js{\ltw,\ltx,\lty}$, has $18$ members.  Each node is labeled by a
phrase---a sum of products.  The dotted equal-sign connects two semicircles,
which represent two instances of the
\emf{same} node, ${\ltw\ltx+\ltw\lty+\ltx\lty}$.  This node is self-dual.
\endgraf
  The lattice exhibits a mirror symmetry across a line passing through nodes
$\ltw$, $\ltx$ and~$\lty$.  Each node in the lattice is mir\-ror-sym\-met\-ric with
its dual.
\jb{The \df{dual} of a phrase is obtained by replacing addition by multiplication,
and vice\,versa.  So the dual of $\ltw+\jp{\ltx\lty}$ is $\ltw\jp{\ltx+\lty}$,
which equals $\ltw\ltx+\ltw\lty$.}
\endrfigure
\par
  Fix an alphabet
$\js{\ltw,\ltx,\lty,\dots,\ltz}$
comprising $\nproto$~letters.  A \df{word} is a non-empty product of letters,
e.g~$\ltw\lty\ltz$.  Each expression built from \logand//\logor/, e.g,
	$$
\alpha \quad\idtb\quad
\jpG{
  \ltx \,+\, \ltw\jpg{\jp{\ltw\ltz+\ltx\ltz}\lty \,+\, \ltx}
}\jpg{\ltw+\lty+\ltz}
\;+\;
\lty\jp{\ltx+\lty}
	 $$%
can be rewritten
\jfo{courtesy of the distributive laws}
as a non-empty sum of words.  Moreover:
\roster
\item"$\bullet$"
  By idempotency and commutativity, each word has no repeated letters.
\item"$\bullet$"
  By absorption, no word is a subword of another word in the sum.
\endroster
Such a reduced sum will be called a \df{phrase}.  The above expression reduces
to the phrase
$\alpha=\ltw\ltx+\ltx\ltz+\lty$.
\par
  Assigning $0$~\jp{={\bf{false}}}~and $1$~\jp{={\bf{true}}} to each of the
$\nproto$ symbols $\ltw,\dots,\ltz$, gives a phrase~$\alpha$ the value
$0$~or~$1$.  So the phrase \jco{thus viewed} is a non-decreasing Boolean
function of its variables.  Conversely, each non-decreasing Boolean function
reduces to a unique phrase.  Consequently:
\quad
\emfs{The free distributive lattice can be written as the lattice of phrases,
or of non-decreasing Boolean functions over $\nproto$~variables.}

\subhead
  \de/ Numbers
\endsubhead
  $\ded\jb\nproto$ is the cardinality of~$\latdede\jb\nproto$.  Some known
values
\jkfootnote{See \cite{Comtet} or {\tt superseeker}.
\quad
  In my definition of $\latdede\jb\nproto$ I omitted two phrases: the
constant~0 function \jp{the empty sum} and the constant~1 function \jp{the sum
whose only term is the empty word}.  Some authors include these phrases, and
so their \de/ numbers are two higher than those listed here.}
are $1, 4, 18, 166, 7579, \allowbreak 7828352, 2414682040996, \allowbreak
56130437228687557907786$.

\par
  A lower bound on the sequence comes from words using half the alphabet.
Let $H\idtb\jf{\nproto/2}$ and
consider those words which use exactly~$H$ of the
$\nproto$~letters.  Evidently every sum of such words is a phrase, and there
are $-1+{\binom\nproto{H}}$ such non-empty sums.
Stirling's approximation to the binomial coefficient ${\binom\nproto{H}}$ gives
	$$
\ded(\nproto)
\quad\ge\quad
{\tfrac12}\cdot 2\!\uparrow\!{\tbinom\nproto{\jf{\nproto/2}}}
\quad\approx\quad
2^{\jbg{2^\nproto \bigm/ \sqrt{\pi \nproto/2}}}
\,.	 $$%
The upshot is that~(\eetable) is a naturally occuring table of numbers which grows
doubly-exponentially in one direction, and {apparently} polynomially in the
other.

\subhead
  Lifting $\minmal\jpg{\tlable\jp\proto}$
\endsubhead
  There is a lattice homomorphism $\phi_1$ from $\latdede=\latdede\jb\nproto$
onto $\prjt\bseta1$.  Simply specify a bijection~$\phi_1$ from the
$\nproto$~generators $\js{\ltw,\dots,\ltz}$ onto the multiset
$\prjt\proto1$, then extend~$\phi_1$ by the two lattice operations.
Similarly, let $\phi_\dira$ be a homomorphism from $\latdede$
onto~$\prjt\bseta\dira$.
\par
  In consequence, the Cartesian product
$\overline\phi\idtb{\phi_1\ct\cdots\ct\phi_\dm}$ is a lattice homomorphism
which \jco{for each direction~$\dira$} makes the following diagram commute.
\group
\def\hold{\latdede\times\cdots\times\latdede}
	$$
    \CD
\hold
  @> {\;\comb_\dira\;} >>
  \hold
    \\
@V {\displaystyle\overline\phi} VV
@V {\displaystyle\overline\phi} VV
    \\
{{\prjt\bseta1}\times\cdots\times{\prjt\bseta\dm}}
  @> {\;\comb_\dira\;} >>
  {{\prjt\bseta1}\times\cdots\times{\prjt\bseta\dm}}
    \endCD
	 $$%
\endgroup
The set of cubes $\brkw\idtb{\ltw^{\ut\dm}}$, \dots,
$\brkz\idtb{\ltz^{\ut\dm}}$ in~${\latdede^{\ut\dm}}$, upstairs, corresponds to
the proto-set~$\proto$, downstairs.
So the homomorphism $\overline\phi$ provides an order-preserving surjection
	$$
\extdfold\jpg{\js{\brkw,\dots,\brkz}}
\;\longrightarrow\;
\extdfold\jpg{\proto}
\,,	 $$%
and consequently each minimal brick in $\extdfold\jp\proto$ comes
from some minimal brick upstairs.
\par
  The upshot is this: Suppose we once-and-for-all compute
$\bsetf\idtb{\bsetm\jpg{\js{\brkw,\dots,\brkz}}}$,
the minimal bricks in
$\ext_{{1\doubledots\dm}}\jpg{\js{\brkw,\dots,\brkz}}$.
Then we know $\bsetm\jp\proto$ for \emf{every} proto-set $\proto$ of
$\nproto$~many $\dm$-bricks:
\;
\emfs{$\bsetm\jp\proto$ comprises the minimal members of the homomorphic image
$\overline\phi\jp\bsetf$.}
\par
  This observation intimates that we might profitably lift our regard to the
product lattice~$\latdede^{\ut\dm}$.

\huchapter3{The Product Lattice}
  We now recast the PC in a lattice setting, stating five Facts
\jco{F1--F5} used for its proof, but leaving their technical demonstration to
the purely algebraic paper~\cite{Kin5}.

\startnormal
  Let $\latprod=\stdprodlat$ denote the $\dm$-fold product
lattice
$\latdede\jb\nproto\times\cdots\times\latdede\jb\nproto$.
Given bricks $\brka,\brkb\in\latprod$, write \dfw{$\brka\cix_\dira\brkb$} for
$\comb_\dira\js{\brka,\brkb}$.  Evidently $\cix_1,\dots,\cix_\dm$ are
CIA~operators, each distributing over every other.
{However}, once $\dm$ exceeds~$2$, the $\cix_\dira$ operators 
\emf{no longer fulfill absorption}.  Here is a $\dm=3$ counterexample.
\group
	$$ \let\holdtimes=\times \def\times{\,\holdtimes\,} \def\;{\,\,\,}
  \align
\brkw \cix_1 \jpg{\brkw \cix_2 \brky}
\;&=\;
\ltw\jp{\ltw{+}\lty} \times \ltw{+}\jp{\ltw\lty} \times \ltw{+}\jp{\ltw{+}\lty}
    \\
\;&=\;
\ltw \times \ltw \times \jp{\ltw{+}\lty}
\;\neq\;
\brkw \,.
    \endalign
	 $$%
\endgroup
Such operators $\js{\cix_\dira}_\dira$ are called \dfw{semilattice operators}
\jco{\cite{Kn\&Ro}}
and the algebraic structure
$\jpg{\latprod, \jp{\cix_1,\dots,\cix_\dm}}$
is a \dfw{multi-semilattice}.

\definition{Definitions}
  The \df{alphabet} of a phrase is the set of letters it uses.  So  the alphabet
of the product $\alpha\idtb\jp{\ltw+\lty\ltz}\jp{\ltw+\ltx+\lty\ltz}$ is
just~$\js{\ltw,\lty,\ltz}$, since $\alpha$ reduces to~$\ltw+\lty\ltz$.
  The alphabet of a {brick}~$\brkb$, written $\alf\jp\brkb$, is the union of
the alphabets of all his sidelengths.  Say that $\brkb$ is \df{balanced} if
all his sidelengths have the same alphabet.
\par
  The sidelengths of $\brkb$ live in a \de/ sublattice of $\latdede\jb\nproto$;
the sublattice generated by~$\alf\jp\brkb$.  The maximum element of this
sublattice, which is the sum of the letters of~$\alf\jp\brkb$, will be called the
\df{envelope} of~$\brkb$.  We write it as~$\env\brkb$.
For example, if $\brkb$ is the $2$-brick
$\jp{\ltx\lty+\ltx\ltz}\times{\jp{\ltz\ltw+\lty}}$, then
$\env\brkb={\ltw+\ltx+\lty+\ltz}$.
\par
   A $\comb$ expression
{which only uses those bricks that are cubes over the given alphabet}, e.g,
	$$
\jp{\brky \cix_{19} \brkw} \cix_{57} \jpg{\brkz \cix_{18} \jp{\brkw \cix_{19} \brkx}}
\tag\put\eegoodexpr
	 $$%
will be called a \df{good} expression.
A \emf{brick}~$\brkb$ is \df{good} if it is the value of some good
expression.  So the set of $\stdprodlat$-good bricks is precisely
	$$
\bsetb_\dm \idtb \ext_{{1\doubledots\dm}}\jpg{\js{\brkw,\dots,\brkz}}
\,.	 $$%
If $\brkb$ is a \emf{minimal} member of $\bsetb_\dm$, say that $\brkb$ is
\df{$\stdprodlat$-minimal}.  Alternatively,
\dfw{$\brkb$ is minimal for~{$\stdprodlat$}}.
\enddefinition

\proclaim{Fact~F1  {\ept\cite{Kin5, Full-Alphabet/Decomp Lemmata}}}\;
  In $\stdprodlat$ suppose brick $\brkb$ is good.  Then $\brkb$ can be built
by some good expression which \emf{only} uses cubes over the alphabet
of~$\brkb$.
\par
  Moreover, if $\brkb$ is minimal then this expression can be chosen to employ
operations $\cix_\dira$ \emf{only} in those directions~$\dira$ where the
sidelength, $\slb_\dira$, is \emff{not} the envelope~$\env\brkb$.
\endproclaim

\proclaim{Fact~F2  {\ept\cite{Kin5, Equal-Alphabet Lemma}}}\;
  If $\brkb$ is $\stdprodlat$-minimal, then $\brkb$ is balanced.
\endproclaim

\let\dimnum=\relax 
\specialhead
  Overall Strategy
\endspecialhead
  An impediment to discussing the function
$\dm\mapsto\maxrank\jp{\nproto,\dm}$ is that the ambient lattice $\stdprodlat$
changes, alas, with~$\dm$.  Further, a natural avenue towards proving PC is by
inducting on~$\nproto$, and here too this inconvenience arises.  In order to
bypass this hindrance, a straightforward approach is to create a big lattice,
$\iilat$, to serve as a common setting for all values of
$\nproto$~and~$\dm$.
\par
  Let ${\seq\ltw\FIR},{\seq\ltw\SEC},\dots$ be an infinite list of letters and let
$\latdede\jb\ii$ be the free distributive lattice that they generate;
$\latdede\jb\ii$ comprises all finite sums of finite words.
Identifying $\latdede\jb\nproto$ with the lattice generated by
${\seq\ltw\FIR},\dots,{\seq\ltw\LAS}$ shows that $\latdede\jb\ii$ is the
direct limit
$\latdede\jb1\hookrightarrow\latdede\jb2\hookrightarrow\dots$ of lattices.
Lastly, for $\nproto$ any value in $\js{1,2,\dots,\ii}$, let $\latprod\jb{\nproto,\ii}$
represent the infinite-product lattice
	$$
\latprod\jb{\nproto,\ii}
\;\idtb\;
\latdede\jb\nproto \times \latdede\jb\nproto \times \cdots
\,.	 $$%

\subhead
  Minimality and Alphabet Size
\endsubhead
  As a first step to implementing the strategy we ask:
\;\emfs{If brick $\brkb$ is $\stdprodlat$-minimal, must he be
minimal for $\latprod\jb{\nprotopo,\dm}$\,?}
\par
  \jq{Yes}, except for the following type of triviality when~$\dm=1$:
The $1$-dimensional brick $\brkw\cix_1\brkx=\ltw\ltx$ is minimal with respect
to alphabet~$\js{\ltw,\ltx}$; it is $\latprod\jb{2,1}$-minimal.  But
$\brkw\cix_1\brkx$ is not $\latprod\jb{3,1}$-minimal, since
$\brkw\cix_1\brkx\cix_1\brky=\ltw\ltx\lty$ is a proper divisor.
\par
  Once $\dm=2$, however, this triviality evaporates.  Now
$\brkw\cix_1\brkx$ equals $\jp{\ltw\ltx}\times\jp{\ltw+\ltx}$, which is
minimal even for $\latprod\jb{\ii,2}$.
\par
  These observations can be interpreted as explaining why the $\dm=\dimnum1$
values of rank $\maxrank\jp{\nproto,\dimnum1}$ do not fit the polynomial
pattern observed in Table~{\eetable}.

\subhead
  Minimality and Dimension
\endsubhead
  Given a good $\stdprodlat$-brick~$\brkb$, there exists some good
expression which fabricates~$\brkb$, using only cubes over~$\alf\jp\brkb$.
This same expression \jco{when interpreted in $\latprod\jb{\nproto,\ii}$}
yields a brick \jco{$\brkb^*$} which is infinite dimensional.  Evidently
	$$
\display
  $\brkb^* \,=\,
\slb_1\times\slb_2\times\cdots\times\slb_\dm
\times {\env\brkb} \times {\env\brkb} \times\cdots$,
where the envelope ${\env\brkb}$ is the sum of the letters in~$\alf\jp\brkb$.
\enddisplay
	 $$%
It turns out that if $\brkb$ is $\stdprodlat$-minimal, then
$\brkb^*$ is $\latprod\jb{\nproto,\ii}$-minimal.
The converse holds
\jkfootnote{For each good brick $\brka$ in $\latprod\jb{\nproto,\ii}$, there is
some integer~$\dm$ such that $\sla_{\dm+1}=\sla_{\dm+2}=\dots$, all being the
envelope of~$\brka$.  Thus each in-fi\-nite-\dimal/ good~$\brka$ is of the
form~$\brkb^*$, for some fi\-nite-\dimal/ brick~$\brkb$.}
once $\dm\ge2$.
\par
  We can summarize these facts as follows.

\proclaim{Fact~F3  {\ept\cite{Kin5,  Universally-minimal Lemma}}}\;
  For each $\dm\ge\dimnum2$ and each~$\nproto$:  Brick~$\brkb$ is
$\stdprodlat$-minimal {\,\;IFF\;\,} $\brkb^*$ is $\iilat$-minimal.
\endproclaim

\par
  Courtesy of~F3, we can henceforth work entirely in the infinite
lattice~$\latprod\idtb\iilat$, and so we should adjust our notation
accordingly.
 Let $\seq\brkw{n}$ be the $\ii$-dimensional
${\seq\ltw{n}}\ct{\seq\ltw{n}}\ct\cdots$
cube.
A brick $\brkb\in\latprod$ is \dfw{\df{\,$\nproto,\dm$-good\,}} if
	$$
\brkb\,\in\,
\ext_{{1\doubledots\dm}}
\jpG{\jsg{
  {\seq\brkw\FIR},{\seq\brkw\SEC}, \dots, {\seq\brkw\LAS}
  }}
\,.	 $$%
Further, $\brkb$ is \dfw{\df{\,$\nproto,\dm$-minimal\,}} if
it is {$\nproto,\dm$-good} \emf{and} is minimal for~$\latprod$.
\extranotes
  \jept{An \jq{$\nproto$-good/minimal} brick is one which is
{$\nproto,\dm$-good/minimal} for some~$\dm$.}
\endextranotes
Redefine
\jkfootnote{This changes the value of $\maxrank\jp{\nproto,\dm}$ only for
$\dm\le1$.}
$\maxrank\jp{\nproto,\dm}$ to now mean the number of
${\nproto,\dm}$-minimal bricks.  For $\nproto=1,2,\dots$, the Polynomial
Conjecture now becomes
	$$
    \gather
{\mskip -20mu}
\text{PC}\jb\nproto:{\mskip10mu}
\text{\sl As $\dm$ takes on the values $0,1,2,\dots$, the resulting
$\maxrank\jp{\nproto,\cdot}$~function,}
    \\
\dm\mapsto\jb{\text{\rm Number of $\nproto,\dm$-minimal bricks}} \; ,
    \\
\foldedtext\foldedwidth{.9\hsize}{\sl
  is a degree-$\jp\nprotomo$ polynomial.%
  }
    \endgather
	 $$%

\subhead
  An Interpretation
\endsubhead
  For $\dm={\dimnum0},{\dimnum1},{\dimnum2}$, the maxrank number
$\maxrank\jp{\nproto,\dm}$ can be regarded as the cardinality of particular
subsets of the \de/ lattice~$\latdede\jb\nproto$.
\group
\roster
\item"$\bullet$"
  $\maxrank\jp{\nproto,\dimnum0}=\nproto$
is the cardinality of the set of \emf{generators} of~$\latdede\jp\nproto$.
\item"$\bullet$"
  Taking the closure of the generating set \jco{under $\cdot$~\jp{product}} 
gives the set of \emf{words}.  Consequently
$\maxrank\jp{\nproto,\dimnum1}={2^\nproto}-1$.
\item"$\bullet$"
  Closing the set of words under~$+$ \jp{sum}, yields~$\latdede\jp\nproto$,
the set of \emf{phrases}.  It is straightforward to check that
	$$
\alpha \;\mapsto\; {\alpha \times {\dual\alpha}}
	 $$%
is a bijection from $\latdede\jp\nproto$ onto $\latprod\jb{\nproto,\dimnum2}$,
where $\dual\alpha$ denotes the dual of~$\alpha$.  Thus
$\maxrank\jp{\nproto,\dimnum2}=\jvg{\latdede\jp\nproto}=\ded\jp\nproto$.
\endroster
\endgroup

\specialhead
  A finite certificate for polynomialness
\endspecialhead
  Now consider a $\iilat$-minimal brick $\brkb$; suppose he can be built with
$\cix_1,\dots,\cix_\dm$.  Let
$\slb_{\dira_1}, \dots,\slb_{\dira_\tdm}$ be an enumeration of the non-envelope
sidelengths of~$\brkb$; {so $\dira_\tdm\le\dm$}.  A non-envelope sidelength
$\slb_\dira$ can be recognized immediately: Since $\brkb$ is balanced,
$\alf\jp{\slb_\dira}=\alf\jp{\env\brkb}$, yet ${\slb_\dira}\neq{\env\brkb}$.
Thus $\slb_\dira$ is not a \dfw{pure sum} of letters ---it must contain a word of
length at least~$2$.
\par
  Because all the $\ext$ operators mutually commute, the permuted brick
	$$
\brkc\idtb
\slb_{\dira_1}\times\slb_{\dira_2}\times\cdots\times\slb_{\dira_\tdm}
\times{\env\brkb}\times{\env\brkb}\times\cdots
	 $$%
is also $\iilat$-minimal.
  Furthermore, courtesy of Fact\,F1, brick~$\brkc$ can be constructed only
using $\cix_1,\dots,\cix_\tdm$.  This number $\tdm$ is what we will call the
\df{true dimension} of~$\brkb$.
For example, the true dimension of~(\eegoodexpr) is three.  Each cube
${\seq\brkw{n}}$ has true dimension zero, and these are the only minimal
bricks with true-dim zero.
\par
  Both bricks $\brkb$ and $\brkc$ are built from the pure sum~$\env\brkb$ and
the multiset
$\js{\slb_{\dira_1}, \ldots, \slb_{\dira_\tdm}}$
of sidelengths.  In order to systematically count the bricks thus-build\-able,
we write this data in a canonical way, by fixing some strict total-order~$\lessdot$ on
$\latdede=\latdede\jb\ii$.
\par
  Letting~$K$ be the number of \emf{distinct} sidelengths in
$\slb_{\dira_1}, \ldots, \slb_{\dira_\tdm}$,
we may rewrite this multiset as
	$$
\sla_1,\ondots{r_1},\sla_1 , \,
\sla_2,\ondots{r_2},\sla_2 ,
\; \cdots \;                ,
\sla_K,\ondots{r_K},\sla_K
\;\, ,	 $$%
where $\sla_i$ is repeated $r_i$~times, the sum $r_1+\dots+r_K$ equals~$\tdm$,
and where $\sla_1\lessdot\sla_2\lessdot\ldots\lessdot\sla_K$.
Thus the expression
$\jbg{e\,;\,
  \sla_1^{\ut r_1} ,
  \sla_2^{\ut r_2} ,
  \dots ,
  \sla_K^{\ut r_K}
}$
tells us what bricks can be build from the multiset.

Conversely,
an expression
$\brka\,\idtb\,\jbg{e\,;\,
  \sla_1^{\ut r_1} ,
  \sla_2^{\ut r_2} ,
  \dots ,
  \sla_K^{\ut r_K}
}$
\jco{formed from members $e,\sla_1,\dots,\sla_K\in\latdede$}
is a \df{archetype} if
\roster
\item"$\bullet$"
  Sidelength $e$~is a pure sum.
\item"$\bullet$"
  Sidelength $\sla_1\lessdot\ldots\lessdot\sla_K$, and none is a pure sum.
\item"$\bullet$"
   $\sla_1^{\ut r_1}\times\cdots\times{\sla_K^{\ut r_K}}\times{e^{\ut\ii}}$
is a minimal brick.
\endroster
Naturally, we call $r_1+\dots+r_K$ the true dimension of~$\brka$ and write
it~$\tdm\jp\brka$.

\subhead
  Counting
\endsubhead
  Now consider a dimension~$\dm$ greater-equal the true-dim
$\tdm=\tdm\jp\brka$.  The number of ways
\jfo{let's call it~$\nbrka\jp\dm$}
of placing $\tdm$~sidelengths together with $\dm-\tdm$~copies of~$e$, into
$\dm$~positions, is expressible by the multinomial coefficient
	$$
\nbrka\jp\dm
\;=\; \binom{\dm}{r_1,\ldots,r_K,\dm-\tdm}
\,\eqnote\,
\frac{\dm!}{{r_1!}\cdots{r_K!}\cdot{\jp{\dm-\tdm}!}}
\,.	 $$%
In consequence, $\nbrka\jp\dm$ is a polynomial in~$\dm$,
	$$
    \aligned
\nbrka\jp\dm
\;&=\;
\binom{\tdm}{r_1,\ldots,r_K}
\cdot \binom{\dm}{\tdm}
     \\
 &
= {\frac{q}{\tdm!}}\;\cdot\;{\dm\cdot\jbg{\dm-1}\cdots\jbg{\dm-\jp{\tdm-1}}}
\;,
    \endaligned
\tag\pf\eecoef
	 $$%
where $q$ is the \emf{integer}~$\binom{\tdm}{r_1,\ldots,r_K}$.  Remark that
this polynomial has degree~$\tdm$, and gives the correct value of
$\nbrka\jp\dm$
\jfo{namely zero}
for each $0\le\dm<\tdm$.
\smallskip
  Say that $\brka$ is an \dfw{$\nproto,\tdm$-archetype}, where $\tdm=\tdm\jp\brka$,
if $\nproto$ is large enough that $\alf\jp\brka$ is a subset of
$\jsg{{\seq\ltw\FIR},\dots,{\seq\ltw\LAS}}$.
Then, with $\nproto$ held fixed,
	$$
\maxrank\jp{\nproto,\dm} \,=\,
\sum_{\tdm=0}^\ii \sum_\brka
\nbrka\jp\dm
\,,	 $$%
where $\brka$ ranges over the \jp{finite} set of $\nproto,\tdm$-archetypes.
Therefore,
the inner sum {$\sum_\brka\nbrka\jp{\cdot}$}
is a polynomial of degree~$\tdm$.  We obtain the following.

\proclaim{Theorem \put\eearchetype}
  The function $\dm\mapsto\maxrank\jp{\nproto,\dm}$ is a polynomial {\,IFF\,}
the set of $\nproto$-archetypes is finite.  In that instance, letting $M$
denote the maximum true dimension taken over the $\nproto$-archetypes, the
polynomial $\maxrank\jp{\nproto,\cdot}$ has degree~$M$.
\endproclaim

\smallskip
  This theorem permits a computer proof of
$\text{PC}\jb\nproto$, for small values of~$\nproto$.  As $\dm=1,2,\dots$,
perform the following:
\group\sl
\roster
\item"({\it i})"
  Having computed $\bsetm_{\dm-1}$
\jco{the set of $\latprod\jb{\nproto,\dm-1}$-minimal bricks}
\dfw{extend} each such brick~$\brkb$ to the $\dm$-brick
{\,}$\brkb'\idtb\brkb\times{\env\brkb}$.
\item"({\it ii})"
  Compute $\bseta$, the set of bricks
$\comb_{\dm}\jpg{\js{\brka',\dots,\brkc'}}$ as $\js{\brka',\dots,\brkc'}$
ranges over all non-void sets of extended bricks.  Thus $\bsetm_\dm$ equals
$\minmal\jp\bseta$.  If some brick in $\bsetm_\dm$ has true-dim~$\dm$, then
GOTO step~{\rm({\it i})}.
\; Otherwise STOP; the maximum true-dimen\-sion of an $\nproto$-archetype
is~$\dm-1$, and collection $\bsetm_{\dm}$ is a certificate of this.
\endroster
\endgroup
  My friend Hugh Redelmeier wrote an intricate computer program to compute
archetypes.  After running for more than a week on the $\nproto=5$~case, his
program constructed all the archetypes and discovered that the maximum
true-dim is~$4$, thus establishing~$\text{PC}\jb{5}$.
The certificate $\bsetm_{5}$ has $273540$~members.

\subhead
  The last ingredient
\endsubhead
  Courtesy of the theorem, the Polynomial Conjecture follows from these two
facts.

\proclaim{Fact F4}
  The brick
	$$
\jpg{\ldots
  \jpg{
    \jp{{\seq\brkw\FIR} \cix_1 {\seq\brkw\SEC}}
    \cix_2 {\seq\brkw\THI}
    }
  \cix_3
  \ldots
  }
\cix_\nprotomo   {\seq\brkw\LAS}
	 $$%
is $\nproto$-minimal.  Furthermore, none of its first $\nprotomo$~sidelengths
is a pure sum, so $\nprotomo$ is indeed the true dimension of the brick.
\endproclaim

\proclaim{Fact F5}
  Each $\nproto$-minimal brick has true dimension at most~$\nprotomo$.
\endproclaim

  This latter result follows from a simultaneous induction on $\nproto$~and
$\dm$ within the $\iilat$~lattice.

\huchapter4{Egress}
   Professors T.\,Hamachi~and Y.\,Tomita recently sent me a preprint
\jco{\cite{Ha\&To}} which develops a new technique to extend the computations
done in \cite{\kbool} for the maxrank numbers in Table~\eetable.  And
\jfo{happily} their results agree with Redelmeier's.
\par
  I warmly thank George Bergman, Kevin Keating, Eric Mendelsohn and Hugh
Redelmeier, as well as the University of Toronto for its hospitality during a
sabbatical visit.

\subhead
  Questions
\endsubhead
  Here are some algebraic questions suggested by the argument.
\par
  Is there a reasonably simple recurrence relation among the
$\maxrank\jp{\nproto,\cdot}$ polynomials?   If so, this would likely lead to a
new method to compute the \de/ numbers, a sequence which has been the object
of considerable study.
\par
  An even more likely place to find a recurrence relation is in the
$2$-parameter table of values~$\text{Arch}\jp{\nproto,\dm}$
\jco{for $\dm<\nproto$}
whose entry is the number of $\nproto,\dm$-archetypes.

\goodbreak
\startnormal
  Affirmative answers to the following would speed up the computation of
archetypes:
\emfs{If minimal bricks $\brkb$ and $\brkc$ have disjoint alphabets, must
$\brkb\cix_1\brkc$ be minimal?  Can each minimal brick $\brkt$ be obtained
\jept{\sl from the given cubes} by a succession of $\cix_\dira$~operations, so
that at every stage the two operand bricks are minimal?}


\define\jct{Jour. Comb. Theory}
\define\amm{Amer. Math. Monthly}

\Refs 
{Mathematical Review numbers \jco{where available} are listed at the end of
each reference.  When {{\sl M\!R:}} numbers are not
available, a call number may be listed.}
\widestnumber\key{dB\&Kla}

\ref
\key{Bar1}
\paper Algebraic theory of brick packing, I
\by F.W.\,Barnes
\jour Discrete Math.
\vol 42
\yr 1982
\pages 7--26
\mathrev{84e:05044a}
  \endref

\ref
\key{Bar2}
\paper Algebraic theory of brick packing, II
\by F.W.\,Barnes
\jour Discrete Math.
\vol 42
\yr 1982
\pages 129--144
\mathrev{84e:05044b}
  \endref

\ref
\key{Comtet}
\book Advanced Combinatorics
\by L.~Comtet
\publ Reidel, Dordrecht, Holland
\yr 1974
\mathrev{57{\,}\#124}
  \endref

\ref
\key {Co\&La}
\paper {Tiling with Polyominoes and Combinatorial Group Theory}
\by J.H.\,Conway \& J.C.\,Lagarias
\jour J. Comb. Theory~A
\yr 1990
\pages 183--208
\vol 53
\mathrev{91a:05030}
  \endref

\ref
\key{Dehn}
\by M.\,Dehn
\paper {\"Uber die Zerlegung von Rechtecken in Rechtecke}
\jour Math. Ann.
\vol 57
\yr 1903
\pages 314--332
  \endref

\ref
\key{Fr\&Ri}
\paper Tiling a Square with Similar Rectangles
\by C.\,Freiling \& D.\,Rinne
\jour Math\. Research Letters
\yr 1994
\pages 547--558
\vol 1
\mathrev{95e:52040}
  \endref

\ref
\key{Gal\&G}
\by D.\,Gale \& C.\,Gardner
\paper On equidecomposability of polygons by rectangles
\finalinfo{Preprint}
  \endref

\ref
\key{Ha\&To}
\by T.\,Hamachi \& Y.\,Tomita
\paper {On the King maximal tiling-rank function}
\finalinfo{Preprint}
  \endref

\ref
\key{Ka\&Sz}
\by G.\,Katona \& D.\,Sz\'asz
\paper {Matching Problems}
\jour \jct
\yr1971
\pages 60-92
\vol 10
\mathrev{46 \#798}
  \endref

\ref
\key{\kk1}
\by K.\,Keating \& J.L.\,King
\paper {Shape Tiling}
\jour Electronic J. Combinatorics	
\yr 1997
\vol 4 {\rm\#2, R12,} 
\pages 48pp.
\mathrev{98e:52023}
  \endref

\ref
\key{\kk2}
\by K.\,Keating \& J.L.\,King
\paper Signed Tilings with Squares
\jour Jour\. of Combinatorial Theory, Series~A
\toappear
  \endref

\ref
\key{\kbool}
\paper {Brick Tiling and Monotone Boolean Functions}
\by J.L.\,King
\finalinfo{Preprint available at: \hfill\allowbreak
 {\tt http://www.math.ufl.edu/$\sim$squash/tilingstuff.html}%
 }
  \endref

\ref
\key Kin2
\paper {Brick Packings and Splittability}
\by J.L.\,King
\finalinfo{Preprint available at: \hfill\allowbreak
 {\tt http://www.math.ufl.edu/$\sim$squash/tilingstuff.html}%
 }
  \endref

\ref
\key Kin5
\paper {Polynomial growth of multi-semilattice rank}
\by J.L.\,King
\finalinfo{In progress}
  \endref

\ref
\key {Kn\&Ro}
\by Arthur Knoebel \& Anna Romanowska
\paper Distributive multisemilattices
\jour Dissertationes Math. (Rozprawy Mat.)
\vol 309
\yr 1991
\pages 42pp.
\mathrev{92g:06008}
  \endref

\ref
\key{Lac\&Sze}
\paper Tilings of the Square with Similar Rectangles
\by M.\,Laczkovich \& G.\,Szekeres
\jour Discrete Comput\. Geom\.
\yr 1995
\pages 569--572
\vol 13
\mathrev{95k:52033}
  \endref

\ref
\key {Propp}
\paper {A pedestrian approach to a method of {C}onway, or, {A} tale of two cities}
\by J.\,Propp
\jour {Mathematics Magazine}
\vol {70{\rm, no\.\,5}}
\yr {1997}
\pages {327--340}
\mathrev{1488869}
  \endref

\ref
\key {Thurst}
\paper Conway's Tiling Groups
\by W.P.\,Thurston
\jour \amm
\yr 1990
\pages 757--773
\vol 97
\mathrev{91k:52028}
  \endref

\endRefs 

\setupfootline

\bye
\enddocument